\documentclass[12pt]{article}

\usepackage{latexsym}
\usepackage{calc}

\usepackage{amssymb}
\usepackage{amsmath}
\usepackage{graphicx}

\usepackage{psfrag}

\usepackage{float}

\usepackage{color}

\usepackage{nicefrac}

\newcounter{hours}\newcounter{minutes}

\def\nr{\par \noindent}
\def\diag{{\rm diag \,}}

\def\beq{\begin{equation}}
\def\eeq{\end{equation}}

\newtheorem{theorem}{Theorem}
\newtheorem{lemma}{Lemma}
\newtheorem{corollary}{Corollary}

\newtheorem{assumption}{Assumption}
\newtheorem{definition}{Definition}

\newtheorem{example}{Example}
\newtheorem{remark}{Remark}
\newcommand{\proof}{\bf Proof: \rm \nr}

\newcommand{\qed}{\hfill $\Box$ \nr \medskip}

\def\ba{\begin{array}}
\def\ea{\end{array}}
\def\beann{\begin{eqnarray*}}
\def\eeann{\end{eqnarray*}}
\def\bea{\begin{eqnarray}}
\def\eea{\end{eqnarray}}

\def\BT{\begin{theorem}}
\def\ET{\end{theorem}}
\def\BL{\begin{lemma}}
\def\EL{\end{lemma}}
\def\BC{\begin{corollary}}
\def\EC{\end{corollary}}
\def\BE{\begin{example}}
\def\EE{\end{example}}
\def\BD{\begin{definition}}
\def\ED{\end{definition}}
\def\BR{\begin{remark}}
\def\ER{\end{remark}}
\def\BAS{\begin{assumption}}
\def\EAS{\end{assumption}}
\def\BI{\begin{itemize}}
\def\EI{\end{itemize}}

\def\BMP{\begin{minipage}{9.5cm}}
\def\EMP{\end{minipage}}
\def\MPT{\begin{minipage}{11.5cm}}
\def\EPT{\end{minipage}}

\def\la{\langle}
\def\ra{\rangle}

\def\R{\mathbb{R}}
\def\E{\mathbb{E}}
\def\P{\mathbb{P}}
\newcommand*\diff{\mathop{}\!\mathrm{d}}
\newcommand*\tr[1]{\mathop{}\!\mathrm{tr}\left(#1\right)}
\renewcommand*\diag[1]{\mathop{}\!\mathrm{diag}\left(#1\right)}

\newcommand*\samethanks[1][\value{footnote}]{\footnotemark[#1]}

\title{Discrete choice prox-functions on the simplex}
\author{D. M\"uller
\thanks{
Department of Mathematics, Chemnitz University of Technology,
Reichenhainer Str. 41, 09126
Chemnitz, Germany; e-mail: david.mueller@mathematik.tu-chemnitz.de, vladimir.shikhman@mathematik.tu-chemnitz.de.
 } \and Yu. Nesterov\thanks{Center for Operations Research and Econometrics (CORE),
Catholic University of Louvain (UCL), 34 voie du Roman
Pays, 1348 Louvain-la-Neuve, Belgium, and National
Research University -- Higher School of Economics, Russia;
e-mail: yurii.nesterov@uclouvain.be. } \and V.
Shikhman\samethanks[1]
}

\begin{document}
\maketitle
\vspace{-5ex}
\abstract{ We derive new prox-functions on the simplex from  
additive random utility models of discrete choice. They are convex conjugates of the corresponding surplus functions. 
In particular, we explicitly derive the convexity parameter of discrete choice prox-functions associated with generalized extreme value models, and specifically with generalized nested logit models.
Incorporated into subgradient schemes, discrete choice prox-functions lead to natural probabilistic interpretations of the iteration steps. As illustration we discuss an
economic application of discrete choice prox-functions in consumer theory. The dual averaging scheme from convex programming naturally adjusts demand within a consumption cycle. 
}

\vspace{2ex}

\section{Introduction}
\label{sec:intro}

The use of prox-functions in convex programming has become standard in recent decades. Originally, they were introduced in the context of mirror descent methods \cite{nemirovski:1983}. 
In order to explain how prox-functions enter into optimization methods, we recall their definition.

\begin{definition}[Prox-function] We say that  $d:\R^n \rightarrow \R\cup\{\infty\}$ is a prox-function on a closed convex set $Q \subset \R^n$ if
\begin{itemize}
    \item[(1)] $d$ is continuous with the domain containing $Q$, i.\,e. $Q \subset \mbox{dom}\, d$.
    \item[(2)] $d$ is strongly convex on $Q$ with respect to a norm $\|\cdot\|$, i.\,e.
    there exists a constant $\beta >0$ such that for all $x, y \in Q$ and $\alpha \in [0,1]$ it holds:
\[
   d(\alpha x +(1-\alpha)y) \leq \alpha d(x) + (1-\alpha) d(y) - \frac{\beta}{2} \alpha (1-\alpha) \|x-y\|^2.
\]
%   We call $\beta$ the convexity parameter of $d$.
    \item[(3)] The computation of the convex conjugate
    \begin{equation}
        \tag{A}\label{eq:pc}
         d^*(s) = \max_{x \in Q} \, \left\langle s, x \right\rangle - d(x)
         \end{equation}
    is simple, i.\,e. the unique maximizer $x(s)$ can be easily obtained for any $s \in \R^n$.
\end{itemize}
\end{definition}
Auxiliary optimization problem (\ref{eq:pc}) 
is known to be a key ingredient for minimizing a convex function $f$ on $Q$ by subgradient schemes. Let us take for $s$ subgradients of $f$ or their weighted aggregates. Then, (\ref{eq:pc}) defines a mapping from the dual space with subgradient information $s$ into the primal space with feasible iterates $x(s) \in Q$. This idea leads in particular to primal-dual subgradient methods for convex problems \cite{nesterov:2013}.

There are at least two advantages of using prox-funcions in (\ref{eq:pc}):
\begin{itemize}
    \item[(a)]  It is well known that the complexity bounds for optimization methods heavily depend on the size of the feasible set $Q$. This value has been traditionally defined with respect to Euclidean norm. However, the size of $Q$, measured with respect to another norm, can be smaller. Thus, by introducing prox-functions, which are strongly convex with respect to an appropriate norm $\|\cdot\|$, it is possible to take into account a particular geometry of the feasible set $Q$.
    \item[(b)] More interestingly, prox-functions often allow natural interpretations of the iteration steps (\ref{eq:pc}) within the convex optimization framework. This feature is important in order to explain agents' behavioral dynamics as being driven by unintentional optimization. 
Let us illustrate this for the well known entropic prox-function 
    \[
       d(p)= \sum_{i=1}^{n} p^{(i)} \ln p^{(i)}
    \]
 on the $(n-1)$-dimensional simplex 
\[
   \Delta=\left\{p \in \R^n \,\left|\, \sum_{i=1}^{n} p^{(i)}  =1, p^{(i)} \geq 0\right.\right\}.
\]
For the feasible set $Q=\Delta$ the auxiliary optimization problem (A) reads:
\[
\max_{p \in \Delta} \, \sum_{i=1}^{n} s^{(i)} p^{(i)} - \sum_{i=1}^{n} p^{(i)} \ln p^{(i)}.
\]
Its unique solution is given by 
\[
    p^{(i)}(s) = \frac{e^{s^{(i)}}}{\displaystyle \sum_{i=1}^{n} e^{s^{(i)}}}, \quad i=1,\ldots,n.
\]
This formula is in accordance with the logit model of discrete choice: $p^{(i)}(s)$ can be viewed as the choice probability of detecting $s^{(i)}$ to be maximal among $s^{(1)}, \ldots, s^{(n)}$.
\end{itemize}

In this paper we introduce prox-functions on the simplex  which allow similar probabilistic interpretation as above. They are derived from the additive random utility models of discrete choice \cite{depalma:1994}. Note that the logit model is just one prominent example within this class. Our main result in Section \ref{sec:prox} states that the convex conjugate of the surplus function, associated with an additive random utility model, is a prox-function on the simplex. For the convex conjugate of the corresponding surplus function we show continuity (Section \ref{ssec:co}), strong convexity (Section \ref{ssec:sc}), and simplicity (Section \ref{ssec:sim}). In particular, we explicitly derive the convexity parameter for the class of
generalized extreme value models \cite{mcfadden:1978}, and specifically of generalized nested logit models \cite{wen:2001}. Section \ref{sec:ea} is devoted to an economic application of discrete choice prox-functions in consumer theory. The discrete choice prox-functions are incorporated into the dual averaging scheme from \cite{nesterov:2013} for consumer's utility maximization. 
 This ensures that the update of internal prices for goods' qualities is due to an additive random utility model.
The dual averaging scheme corresponds to a natural consumption cycle which successively leads to an optimal consumption of goods (Section \ref{sec:cyc}). We mention that the proposed consumption cycle generalizes \cite{nesterov:2016} where the update of internal prices is due to the logit model.

{\bf Notation.} Our notation is quite standard. We denote by $\R^n$ the space of $n$-dimensional column vectors $x =
\left(x^{(1)}, \dots , x^{(n)}\right)^T$, by $\R^n_+$ the set of all vectors with nonnegative components. 
If the components of $x \in \R^n$ are nonnegative (positive), we write $x \geq 0$ ($x > 0$). For $x \in \R^n$ we write $x^{(-i)} \in \R^{n-1}$ meaning that the $i$-th component of $x$ is missing. Analogously, we write $x^{(-i,j)} \in \R^{n-2}$ meaning that both  $i$-th and $j$-th components of $x$ are missing.
For $x, y \in \R^n$ we introduce the standard scalar product and 
(if additionally $y > 0$) the vector division:
\[
\la x, y \ra = \sum\limits_{i=1}^n x^{(i)} y^{(i)}, \quad
\frac{x}{y} = \left(\frac{x^{(1)}}{y^{(1)}}, \ldots, \frac{x^{(n)}}{y^{(n)}}\right)^T.
\]
For $x \in \R^n$ we use the following norms:
\[
\|x\|_1 = \sum_{i=1}^{n} \left| x^{(i)}\right|, \quad \|x\|_\infty= \max_{1 \leq i \leq n} \left| x^{(i)}\right|.
\]
Note that they are dual to each other, i.\,e.
\[
  \|x\|_\infty = \sup_{\left\|y \right\|_1 \leq 1} \langle y,x\rangle, \quad
  \|x\|_1 = \sup_{\left\|y \right\|_\infty \leq 1} \langle x,y\rangle.
\]
We denote by $e_j \in \R^n$ the $j$-th coordinate vector of $\R^n$. All components of the vector $e \in \R^n$ are equal to one.
%The same notation is used for different spaces, which are always determined by the context.
The space of $(n \times n)$-matrices with real-valued entries is denoted by $\R^{n \times n}$. We use the induced matrix norm for $A \in \R^{n\times n}$: 
\[
\left\| A\right\|_{\infty,1} = \max_{\|z\|_\infty \leq 1} \left\| A z\right\|_1.
\]
Given a twice differentiable function $f:\R^n \rightarrow \R$, $\nabla f$ denotes the gradient and $\nabla^2 f$ stands for the Hessian matrix.

\section{Discrete choice prox-functions on the simplex}
\label{sec:prox}

We derive discrete choice prox-functions on the simplex from additive random utility models.

\subsection{Additive random utility models}
\label{sec:aru}

The additive random utility framework has been first introduced in economic context \cite{mcfadden:1978}. It aims to model the discrete choice from a finite number of alternatives $\{1, \ldots, n\}$ by a rational decision-maker prone to some random errors. Accordingly, the $i$-th alternative is endowed with the utility  
\[
   u^{(i)} + \epsilon^{(i)},
\]
where $u^{(i)} \in \mathbb{R}$ is its deterministic part and $\epsilon^{(i)}$ is a random error. We denote by
\[
u = \left(u^{(1)}, \ldots, u^{(n)}\right)^T, \quad \epsilon = \left(\epsilon^{(1)}, \ldots, \epsilon^{(n)}\right)^T
\]
the vectors of deterministic utilities and of random utility shocks, respectively.  
The following assumption on the stochastic errors is standard, see e.\,g. \cite{depalma:1994}.

\begin{assumption}
\label{ass:rnd}
The random vector $\epsilon$ follows a joint distribution with finite mean that is absolutely continuous with respect to the Lebesgue measure and fully supported on $\R^n$.
\end{assumption} 
Since a rational decision-maker chooses alternatives with the maximal utility, the corresponding surplus is given by the expectation
\[
  E(u) = \E_\epsilon \left(\max_{1 \leq i \leq n} u^{(i)} + \epsilon^{(i)} \right).
\]
It is well-known that the surplus function $E$ is convex and differentiable \cite{depalma:1994}. In particular, its partial derivatives can be expressed as choice probabilities:
\begin{equation}
 \label{eq:der1}
  \frac{\partial E(u)}{\partial u^{(i)}} = \P \left( u^{(i)} + \epsilon^{(i)} = \max_{1 \leq i \leq n} u^{(i)} + \epsilon^{(i)}\right), \quad i=1, \ldots, n.  
\end{equation}
The latter means that the $i$-th partial derivative of $E$ corresponds to the probability of perceiving the $i$-th alternative as one with the maximal utility among the others. This result is known  as the Williams-Daly-Zachary theorem in the discrete choice literature \cite{mcfadden:1978, mcfadden:1981}. The formula (\ref{eq:der1}) is valid due to the fact that, under Assumption \ref{ass:rnd}, the ties between the alternatives occur with zero-probability, i.\,e.
\[
   \P\left(\epsilon^{(i)} - \epsilon^{(j)} = c \right) =0 \quad 
   \mbox{for all } i \not = j \mbox{ and } c \in \R.
\]

\subsection{Convex conjugate of the surplus function}
\label{sec:cc}
We turn our attention to the convex conjugate $E^*:\R^n \rightarrow \R\cup\{\infty\}$ of the surplus function:
\[
   E^*(p) = \sup_{ u \in \R^n} \, \langle p, u \rangle - E(u),
\]
where $p=\left(p^{(1)}, \ldots, p^{(n)}\right)^T \in \R^n$ is the vector of dual variables.

\subsubsection{Continuity}
\label{ssec:co}
We discuss the continuity of the convex conjugate $E^*$ on its domain
\[
  \mbox{dom}\, E^* = \left\{ p \in \R^n \,|\, E^*(p) < \infty \right\}.
\]
For that, we need some elementary properties of the surplus function $E$ listed below.

\begin{lemma}[Elementary properties of $E$] 
\label{lem:el}
For the surplus function $E$ it holds:
\begin{itemize}
    \item[(E1)] $E(u+\gamma e) = E(u) + \gamma$ for all $\gamma \in \R, u \in \R^n$.
    \item[(E2)] $E(u) \geq E(v)$ for all $u, v \in \R^n$ with $u \geq v$.
    \item[(E3)] $E(u) \geq \displaystyle \max_{1 \leq i \leq n} u^{(i)} + \min_{1 \leq i \leq n} \E_\epsilon \left(\epsilon^{(i)} \right)$ for all $u \in \R^n$.
\end{itemize}
\end{lemma}
\proof 
\begin{itemize}
    \item[(E1)] The linearity of the expectation provides:
    \[
    E(u+\gamma e) = \E_\epsilon \left(\max_{1 \leq i \leq n} \left(u^{(i)} + \gamma + \epsilon^{(i)}\right) \right) = 
    \E_\epsilon \left(\max_{1 \leq i \leq n} \left( u^{(i)} + \epsilon^{(i)} \right) + \gamma \right)=
    E(u) + \gamma. 
    \]
    \item[(E2)] The monotonicity of the expectation provides:
    \[
     E(u) =  \E_\epsilon \left(\max_{1 \leq i \leq n} u^{(i)} + \epsilon^{(i)} \right) \geq 
      \E_\epsilon \left(\max_{1 \leq i \leq n} v^{(i)} + \epsilon^{(i)} \right) = E(v).
     \]
    \item[(E3)] Due to the finite mean condition from Assumption \ref{ass:rnd}, we have for every $i \in \{1, \ldots, n\}$:
    \[
    E(u) = \E_\epsilon \left(\max_{1 \leq i \leq n} u^{(i)} + \epsilon^{(i)} \right) \geq 
    \E_\epsilon \left(u^{(i)} + \epsilon^{(i)} \right)  \geq u^{(i)} + \min_{1 \leq i \leq n}\E_\epsilon \left(\epsilon^{(i)} \right).
    \]
\end{itemize}
\qed

\begin{theorem}[Continuity of $E^*$]
\label{th:ct}
   The convex conjugate $E^*$ is continuous on its domain $\mbox{dom}\, E^*$ which coincides with the simplex $\Delta$.
\end{theorem}
   
\proof Let us first show that $\mbox{dom}\, E^* \subseteq \Delta$. For $p \in \R^n$ with $\langle p,e\rangle \not =1$ we have:
\[
  E^*(p) \geq \sup_{ \gamma \in \R}\, \langle p,v + \gamma e \rangle - E(v + \gamma e) \overset{\mbox{(E1)}}{=}  \langle p,v \rangle - E (v) + \sup_{ \gamma \in \R} \, \gamma (\langle p,e \rangle - 1) = \infty,
\]
where $v \in \R^n$ is fixed. For $p \in \R^n$ with $p^{(i)} < 0$ for an $i \in \{1, \ldots, n\}$ we have:
\[
 E^*(p) \geq \sup_{ \gamma \leq 0} \, \langle p, \gamma e_i \rangle - E\left(\gamma e_i\right) \overset{\mbox{(E2)}}{\geq}
 \sup_{ \gamma \leq 0} \, \gamma p^{(i)} - E\left(0\right) = \infty,
\]
where $e_i$ denotes the $i$-th coordinate vector. 
Secondly, we prove that $\mbox{dom}\, E^* \supseteq \Delta$. For that, it is sufficient to show that $E^*$ is bounded from above on $\Delta$. Due to (E3) from Lemma \ref{lem:el}, it holds:
\begin{equation}
\label{eq:ub}
\begin{array}{rcl}
\displaystyle \sup_{p \in \Delta} \, E^*(p) &=& \displaystyle \sup_{p \in \Delta} \left(\sup_{ u \in \R^n} \, \langle p, u \rangle - E(u) \right) =
  \sup_{ u \in \R^n} \left( \sup_{p \in \Delta} \, \langle p, u \rangle - E(u) \right) \\ \\
&=& \displaystyle 
  \sup_{ u \in \R^n} \left(\max_{1 \leq i \leq n} u^{(i)} - E(u)\right) \leq - \min_{1 \leq i \leq n} \E_\epsilon \left(\epsilon^{(i)} \right).
\end{array}
\end{equation}
Further, we discuss the continuity of $E^*$ on the simplex $\Delta$. 
Since $E^*$ is convex, it is continuous on the relative interior $\mbox{rint}(\Delta)$ of its domain. The continuity of $E^*$ on the whole domain $\Delta$ can be deduced by an application of the Gale-Klee-Rockafellar theorem.
The Gale-Klee-Rockafellar theorem says that a convex function is upper semi-continuous at every point at
which its domain is polyhedral \cite{gkr:1968}. Note that the domain of $E^*$ -- the $(n-1)$-dimensional simplex $\Delta$ -- is polyhedral. Moreover, the convex conjugate of a function is always lower semi-continuous, so is $E^*$ on $\Delta$. Together, the lower and upper semi-continuity of $E^*$ on $\Delta$ provides the claim. \qed

Theorem \ref{th:ct} says that the convex conjugate $E^*$ is finite on the simplex $\Delta$.
The latter can be viewed as the set of probability distributions. Hence, the dual variables $p$ can be interpreted as the probabilities attached to the alternatives $\{1, \ldots, n\}$.

\begin{corollary}[Upper bound for $E^*$] 
\label{cor:ub}
The convex conjugate $E^*$ is bounded from above on its domain $\Delta$, namely it holds:
\[ 
   E^*(p) \leq - \min_{1 \leq i \leq n} \E_\epsilon \left(\epsilon^{(i)} \right) \quad \mbox{for all } p \in \Delta.
\]
\end{corollary}

\proof The assertion follows from the derivation in (\ref{eq:ub}). \qed

\subsubsection{Strong Convexity}
\label{ssec:sc}

We show that the convex conjugate $E^*$ is strongly convex under suitable assumptions, and estimate its convexity parameter.

\begin{definition}[Strong convexity of $E^*$] The convex conjugate $E^*: \Delta \rightarrow \R$ is $\beta$-strongly convex with respect to the $\|\cdot\|_1$ norm if for all $p, q \in \Delta$ and $\alpha \in [0,1]$ we have:
\[
  E^*(\alpha p +(1-\alpha)q) \leq \alpha E^*(p) + (1-\alpha) E^*(q) - \frac{\beta}{2} \alpha (1-\alpha) \|p-q\|_1^2.
\]
The positive constant $\beta$ is called the convexity parameter of $E^*$.
\end{definition}

The strong convexity of $E^*$ is closely related to the strong smoothness of $E$.

\begin{definition}[Strong smoothness of $E$] The surplus function $E: \R^n \rightarrow \R$ is $L$-strongly smooth with respect to the maximum norm $\|\cdot\|_\infty$ if for all $u, v \in \R^n$ we have:
\[
  E(u+v) \leq E(u) + \langle \nabla E(u), v \rangle + \frac{L}{2} \|v\|_\infty^2.
\]
The positive constant $L$ is called the smoothness parameter of $E$.
\end{definition}

The following duality result between the strong convexity of $E^*$ and the strong smoothness of $E$ can be easily deduced from \cite[Theorem 6]{kak:2009} shown there in the general setting.

\begin{lemma}[Strong convex/smooth duality]
\label{lem:sd}
The convex conjugate $E^*$ is $\beta$-strongly convex with respect to the $\|\cdot\|_1$ norm if and only if 
the surplus function $E$ is $\frac{1}{\beta}$-strongly smooth with respect to the maximum norm $\|\cdot\|_\infty$.
\end{lemma}

\proof We apply \cite[Theorem 6]{kak:2009} which says that a closed and convex function is $\beta$-strongly convex with respect to a norm if and only if its convex conjugate is $\frac{1}{\beta}$-strongly smooth with respect to the dual norm. For that, we note that $E^*$ is proper and lower semi-continuous, hence, closed. Moreover, by the Fenchel-Moreau theorem we have
\[
   E^{**} = E,
\]
since $E$ is, in particular, a proper, lower semi-continuous, and convex function. Finally, the dual of the $\|\cdot\|_1$ norm is the maximum norm $\|\cdot\|_\infty$.
\qed

In view of Lemma \ref{lem:sd}, we may focus on the strong smoothness of $E$. For the characterization of the latter property, we use the fact that the surplus function $E$ is twice differentiable. Let us compute the second order partial derivatives of $E$. Recall that its $i$-th partial derivative can be written as the choice probability
\[
 \begin{array}{rcl}
        \displaystyle \frac{\partial E(u)}{\partial u^{(i)}} &=& \displaystyle \P \left(\epsilon^{(-i)} - \epsilon^{(i)} \leq u^{(i)} -u^{(-i)} \right) \\ \\ &=& \displaystyle \int_{-\infty}^{u^{(i)} - u^{(-i)}} \int_{-\infty}^{\infty} f_\epsilon \left(y^{(-i)} +x^{(i)}, x^{(i)}\right) \diff x^{(i)} \diff y^{(-i)}, \end{array}
\]
where $f_\epsilon$ is the probability density function of the random utility shocks $\epsilon$. Here, the inner integral is the probability density function of the $(n-1)$-dimensional vector of random differences $\epsilon^{(-i)}-\epsilon^{(i)}$ \cite{depalma:1994}. By differentiating this formula with respect to $u^{(j)}$ for $j\not =i$, we obtain mixed partial derivatives of $E$:
\[
   \frac{\partial^2 E(u)}{\partial u^{(i)} \partial u^{(j)}} = - 
   \int_{-\infty}^{u^{(i)} - u^{(-i,j)}} \int_{-\infty}^{\infty} f_\epsilon \left(y^{(-i,j)} +x^{(i)}, u^{(i)}-u^{(j)}+x^{(i)}, x^{(i)}\right) \diff x^{(i)} \diff y^{(-i,j)}.
\] 
This integral can be interpreted as the probability density that $\epsilon^{(j)} - \epsilon^{(i)}=u^{(i)}-u^{(j)}$, and $\epsilon_{-i,j} - \epsilon^{(i)} \leq u^{(i)}-u^{(-i,j)}$, i.\,e. both alternatives $i$ and $j$ yield the maximal utility.
Analogously, we obtain the second order partial derivative of $E$ with respect to $u^{(i)}$:
\[
   \frac{\partial^2 E(u)}{\partial u^{{(i)}2}} = \sum_{j \not =i}
\int_{-\infty}^{u^{(i)} - u^{(-i,j)}} \int_{-\infty}^{\infty} f_\epsilon \left(y^{(-i,j)} +x^{(i)}, u^{(i)}-u^{(j)}+x^{(i)}, x^{(i)}\right) \diff x^{(i)} \diff y^{(-i,j)}.
\]

\begin{lemma}[$C^2$-characterization of strong smoothness of $E$] 
\label{lem:c2} The surplus function $E$ is $L$-strongly smooth with respect to the maximum norm $\|\cdot\|_\infty$ if for all $u \in \R^n$ it holds:
\[
  \left\| \nabla^2 E(u)\right\|_{\infty,1} \leq L.
\]

\proof For any $u, v \in \R^n$ we have:
\[ 
\nabla E(u) - \nabla E(v) = \int_0^1 \diff \nabla E(v + \tau (u-v)) =
\int_0^1 \nabla^2 E(v + \tau (u-v)) \cdot (u-v) \diff{\tau}.
\]
Hence, the gradient of $E$ is Lipschitz continuous:
 \[
 \begin{array}{rcl}
        \left\|\nabla E(u) - \nabla E(v)\right\|_1 &\leq& \displaystyle
   \int_0^1 \left\| \nabla^2 E(v + \tau (u-v))\cdot(u-v) \right\|_1 \diff{\tau}  \\ \\ 
   &\leq& \displaystyle
   \int_0^1 \left\| \nabla^2 E(v + \tau (u-v))\right\|_{\infty,1} \cdot \left\|u-v\right\|_\infty \diff{\tau} \leq   L \cdot \left\|u-v\right\|_\infty.
\end{array}
    \]
Further, we have:
\[
  E(u+v) - E(u) = \int_0^1 \diff E(u+\tau v) =
  \int_0^1 \langle \nabla E(u+\tau v), v\rangle \diff \tau.
\]
Due to the Lipschitz continuity of $\nabla E$, we obtain:
\[
\begin{array}{rcl}
  E(u+v) - E(u) - \langle \nabla E(u), v\rangle &=& \displaystyle \int_0^1 \langle \nabla E(u+\tau v) - \nabla E(u), v\rangle \diff \tau \\ \\ 
   &\leq& \displaystyle \int_0^1 \left\| \nabla E(u+\tau v) - \nabla E(u) \right\|_1 \cdot \left\|v\right\|_\infty \diff \tau \\ \\ 
   &\leq& \displaystyle
   \int_0^1 L \cdot \left\| u+\tau v - u \right\|_\infty \cdot \left\|v\right\|_\infty \diff \tau = \frac{L}{2} \left\|v\right\|^2_\infty.
  \end{array}
\]
\qed
\end{lemma}

Now, let us consider the set $\mathcal{A}$ of symmetric matrices $A=\left(a_{ij}\right) \in \R^{n\times n}$ satisfying:
\begin{itemize} 
    \item[(A1)] $\displaystyle a_{ii} \geq 0$ for all $i = 1, \ldots, n$, and $a_{ij} \leq 0$ for all $i \not = j$,
    \item[(A2)] $\displaystyle a_{ii} + \sum_{j\not = i}^{} a_{ij} = 0$ for all $i = 1, \ldots, n$.
\end{itemize}
Note that the set $\mathcal{A}$ is closed under matrix addition and multiplication by nonnengative scalars, i\,e. for any $A,B \in \mathcal{A}$ and $\lambda \geq 0$ it holds:
\[
   A+B, \lambda A \in \mathcal{A}.
\]
Moreover, the Hessian matrix $\nabla^2 E(u)$ of the surplus function is an element of $\mathcal{A}$.

\begin{lemma}[Representation of $\|\cdot\|_{\infty,1}$]
\label{lem:nr}
For $A \in \mathcal{A}$ it holds:
\begin{equation}
    \|A\|_{\infty,1} = 4 \max \, \left\{ \left.\sum_{i,j \in K} a_{ij} \,\right|\, K \subset \{1, \ldots, n\}, |K| \leq \left\lfloor \frac{n}{2} \right\rfloor  \right\}.
\end{equation}

\proof The maximum in the definition of the matrix norm 
\[
  \|A\|_{\infty,1} = \max_{\|z\|_\infty \leq 1} \left\| A z\right\|_1 
\]
is attained at some vertex of the feasible set $\left\{z \in \R^n \,\left|\, \|z\|_\infty \leq 1 \right. \right\}$. These are the vectors
\[
   z_K = e_K - e_{K^c},
\]
where $e_K \in \R^n$ is the indicator vector of the subset $K\subset \{1, \ldots, n\}$, i.\,e. $e_K^{(i)}=1$ if $i \in K$, and $e_K^{(i)}=0$ if $i \not \in K$.
We may restrict the choice of $K$ by the condition $|K| \leq \left\lfloor \frac{n}{2} \right\rfloor$, since
\[
  \left\| A z_{K^c}\right\|_1 = 
  \left\| A \left(-z_{K}\right)\right\|_1 = 
\left\| A z_{K}\right\|_1.
\]
For such a fixed subset $K$ we compute
\[
\begin{array}{rcl}
 \displaystyle \left\| A z_{K}\right\|_1 &=& \displaystyle \sum_{i \in K} \left| a_{ii} + \sum_{j \in K\backslash \{i\}} a_{ij} - \sum_{j \in K^c} a_{ij} \right|+  \sum_{i \in K^c} \left| \sum_{j \in K} a_{ij} - a_{ii} - \sum_{j \in K^c\backslash \{i\}} a_{ij} \right| \\ \\ 
   &\overset{\mbox{(A2)}}{=}& \displaystyle 2 \sum_{i \in K} \left| a_{ii} + \sum_{j \in K\backslash \{i\}} a_{ij} \right|+ 2 \sum_{i \in K^c} \left| \sum_{j \in K} a_{ij} \right| \\ \\ 
   &\overset{\mbox{(A1)}}{=}& \displaystyle 2 \sum_{i \in K} \left(a_{ii} + \sum_{j \in K\backslash \{i\}} a_{ij} \right)- 2 \sum_{i \in K^c} \sum_{j \in K} a_{ij} \\ \\ 
   &=& \displaystyle
   2 \sum_{j \in K} \left(a_{jj} + \sum_{i \in K\backslash \{j\}} a_{ij} - \sum_{i \in K^c}  a_{ij} \right) \\ \\ 
   &\overset{\mbox{(A2)}}{=}& \displaystyle
   2 \sum_{j \in K} \left(2 a_{jj} + 2 \sum_{i \in K\backslash \{j\}} a_{ij} \right) = 
   %4 \left\langle e^K, A e^K \right \rangle
   4 \sum_{i,j \in K} a_{ij}.
  \end{array}
\]  
\qed
\end{lemma}

We use Lemma \ref{lem:nr} to estimate the $\|\cdot\|_{\infty,1}$ norm on a particular subset of $A$, which will appear in what follows. 

\begin{corollary} 
\label{cor:pp}
For $p\in \Delta$ we define the matrix $R=\diag{p} - p\cdot p^T$. It holds for the latter:
\[
R \in \mathcal{A}, \quad \|R\|_{\infty,1} \leq 1.
\]
\proof The symmetric matrix $R=(r_{ij})$ fulfills (A1), since in view of $p \in \Delta$ it holds for all $i = 1, \ldots, n$, and $i \not = j$:
\[
   r_{ii}=p^{(i)}\left(1-p^{(i)}\right) \geq 0, \quad r_{ij}=-p^{(i)}p^{(j)} \leq 0.
\]
It also fulfills (A2) due to the following derivation:
\[
   R \cdot e = \diag{p}\cdot e - p\cdot p^T \cdot e = p - p= 0.
\]
Let us fix a subset of indices $K \subset \left\{1, \ldots,n \right\}$.
We estimate the following expression from Lemma \ref{lem:nr} uniformly for all $p \in \Delta$:
\[
    \sum_{i,j \in K} r_{ij} = \sum_{i \in K} p^{(i)} \left( 1- \sum_{j \in K} p^{(j)} \right).
\]
For that, let us solve the maximization problem
\[
    \max \, \sum_{i \in K} p^{(i)} \left( 1- \sum_{i \in K} p^{(i)} \right) \quad \mbox{s. t.} \quad \sum_{i \in K} p^{(i)} \leq 1, \quad p^{(i)} \geq 0 \mbox{ for all } i \in K.
\]
Without loss of generality, we may assume that for its solution holds: 
\[
  \sum_{i \in K} p^{(i)} < 1, \quad p^{(i)} > 0 \mbox{ for all } i \in K.
\]
Otherwise, the optimal value vanishes or we pass over to a smaller subset of indices. The first order optimality condition reads:
\[
  1- \sum_{i \in K} p^{(i)} - \sum_{i \in K} p^{(i)} =0.
\]
Hence, we get $\displaystyle \sum_{i \in K} p^{(i)} = \frac{1}{2}$,
and the optimal value is
\[
  \sum_{i \in K} p^{(i)} \left( 1- \sum_{i \in K} p^{(i)} \right) = 
  \sum_{i \in K} p^{(i)} \left( 1- \frac{1}{2} \right) = \frac{1}{4}.
\]
The application of Lemma \ref{lem:nr} provides the assertion. \qed
\end{corollary}

\begin{lemma}[Estimation of $\|\cdot\|_{\infty,1}$] 
\label{lem:tr}
For $A \in \mathcal{A}$ it holds:
\[
  \|A\|_{\infty,1} \leq 2 \tr{A},
\]
where $\tr{A}$ denotes the trace of the matrix $A$.

\proof For any $z\in \R^n$ it holds:
\[
  \left\| A z\right\|_1 = \sum_{i=1}^{n} \left| \sum_{j=1}^{n} a_{ij} z^{(j)} \right|
\leq \sum_{i=1}^{n} \sum_{j=1}^{n} \left|a_{ij}\right| \cdot \left|z^{(j)}\right| 
\leq \sum_{i=1}^{n} \sum_{j=1}^{n} \left|a_{ij}\right| \cdot \left\|z\right\|_\infty. 
\]
Additionally, we have:
\[
 \sum_{i=1}^{n} \sum_{j=1}^{n} \left|a_{ij}\right|
 \overset{\mbox{(A1)}}{=} \sum_{i=1}^{n} \left( a_{ii} - \sum_{j\not = i}^{} a_{ij} \right) \overset{\mbox{(A2)}}{=} \sum_{i=1}^{n} \left( a_{ii} + a_{ii} \right) = 2 \tr{A}.
\]
Overall, we obtain the inequality:
\[
  \|A\|_{\infty,1} = \max_{\|z\|_\infty \leq 1} \left\| A z\right\|_1 \leq
  \max_{\|z\|_\infty \leq 1} 2 \tr{A} \left\|z\right\|_\infty = 2 \tr{A}.
\]
\qed
\end{lemma}

Now, we are ready to state the general result on the strong convexity of $E^*$. It is given in terms of the differences $\epsilon^{(j)} - \epsilon^{(i)}$, $i \not = j$, of random utility shocks. We recall that the density function of $\epsilon^{(j)} - \epsilon^{(i)}$ can be written as
\[
  g_{i,j} \left(z\right)= \int_{-\infty}^{\infty} \int_{-\infty}^{\infty} f_\epsilon \left(y^{(-i,j)} +x^{(i)}, z +x^{(i)},x^{(i)}\right) \diff x^{(i)} \diff y^{(-i,j)}.
\]
Any point $\bar z_{i,j} \in \R$ which maximizes the density function $g_{i,j}$ is called a mode of the random variable $\epsilon^{(j)} - \epsilon^{(i)}$.

\begin{theorem}[Strong convexity of $E^*$]
\label{th:sc}
   Let the differences $\epsilon^{(j)} - \epsilon^{(i)}$  of random utility shocks have modes $\bar z_{i,j} \in \R$, $i \not = j$.
   Then, the corresponding convex conjugate $E^*$ is $\beta$-strongly convex with respect to the $\|\cdot\|_1$ norm, where the convexity parameter is given by
\[
  \beta = \frac{1}{\displaystyle 2\sum_{i=1}^{n} \sum_{j\not =i} g_{i,j} \left(\bar z_{i,j}\right)}.
\]
\proof 
We estimate the second order derivative of the surplus function $E$ with respect to $u^{(i)}$:
\[
\begin{array}{rcl}
 \displaystyle \frac{\partial^2 E(u)}{\partial u^{{(i)}2}} &=& \displaystyle \sum_{j \not =i} \int_{-\infty}^{u^{(i)} - u^{(-i,j)}} \int_{-\infty}^{\infty} f_\epsilon \left(y^{(-i,j)} +x^{(i)}, u^{(i)}-u^{(j)}+x^{(i)}, x^{(i)}\right) \diff x^{(i)} \diff y^{(-i,j)}\\ \\ 
   &\leq& \displaystyle \sum_{j \not =i} \int_{-\infty}^{\infty} \int_{-\infty}^{\infty} f_\epsilon \left(y^{(-i,j)} +x^{(i)}, u^{(i)}-u^{(j)}+x^{(i)}, x^{(i)}\right) \diff x^{(i)} \diff y^{(-i,j)} \\ \\
   &=& \displaystyle \sum_{j \not =i}  g_{i,j} \left(u^{(i)}-u^{(j)}\right) \leq 
   \sum_{j \not =i}  g_{i,j} \left(\bar z_{i,j}\right).
  \end{array}
\]
Hence, we obtain for the trace of $\nabla^2 E$ the following inequality:
\[
  \tr{\nabla^2 E(u)} = \sum_{i=1}^{n} \frac{\partial^2 E(u)}{\partial u^{{(i)}2}} \leq 
  \sum_{i=1}^{n} \sum_{j\not =i} g_{i,j} \left(\bar z_{i,j}\right).
\]
In view of Lemma \ref{lem:tr}, whose application is justified by the fact that the matrix $\nabla^2 E(u)$ fulfills (A1)--(A2), it holds:
\[
  \left\|\nabla^2 E(u)\right\|_{\infty,1} \leq 2 \tr{\nabla^2 E(u)} \leq 2 \sum_{i=1}^{n} \sum_{j\not =i} g_{i,j} \left(\bar z_{i,j}\right) = \frac{1}{\beta}.
\]
The latter provides, due to Lemma \ref{lem:c2}, that the surplus function $E$ is $\frac{1}{\beta}$-strongly smooth with respect to the maximum norm $\|\cdot\|_\infty$. Finally, we apply Lemma \ref{lem:sd} to conclude that the convex conjugate $E^*$ is $\beta$-strongly convex with respect to the $\|\cdot\|_1$ norm.\qed
\end{theorem}

\begin{remark}[Existence of modes] \textup{We note that the condition on the existence of modes in Theorem \ref{th:sc} cannot be weakened. This can be seen already in case of two alternatives, i.\,e. $n=2$. Then, the second order derivative of the surplus function $E$ is
\[
\nabla^2 E(u) = \left( \begin{array}{rr}
    g_{1,2}\left(u^{(1)}-u^{(2)}\right) & -g_{1,2}\left(u^{(1)}-u^{(2)}\right) \\
    -g_{2,1}\left(u^{(2)}-u^{(1)}\right) & g_{2,1}\left(u^{(2)}-u^{(1)}\right)
\end{array}\right).
\]
After a moment of reflection we realize that  
\[
  g_{1,2}\left(u^{(1)}-u^{(2)}\right) = g_{2,1}\left(u^{(2)}-u^{(1)}\right).
\]
Due to Lemma \ref{lem:nr}, it holds:
\[
  \left\|\nabla^2 E(u)\right\|_{\infty,1} = 4 g_{1,2}\left(u^{(1)}-u^{(2)}\right).
\]
From now on we assume that $g_{1,2}$ is continuous.
Hence, the reverse implication in Lemma \ref{lem:c2} becomes valid, and $E$ is strongly smooth -- or, equivalently, $E^*$ is strongly convex -- if and only if the density function $g_{1,2}$ is bounded on $\R$. The latter property can be characterized by the existence of a mode $\bar z_{1,2}$ of $\epsilon^{(2)}-\epsilon^{(1)}$. Additionally, the convexity parameter of $E^*$ from Theorem \ref{th:sc} can be expressed as 
\[
  \beta = \frac{1}{\displaystyle 2\sum_{i=1}^{n} \sum_{j\not =i} g_{i,j} \left(\bar z_{i,j}\right)} = \frac{1}{\displaystyle 4  g_{1,2}\left(\bar z_{1,2}\right)} = \frac{1}{\left\|\nabla^2 E\left(\bar u\right)\right\|_{\infty,1}},
\]
where the utilities $\bar u = \left(\bar u^{(1)},\bar u^{(2)} \right)^T$ are chosen to satisfy $\bar u^{(1)}-\bar u^{(2)}=\bar z_{1,2}$. This formula provides that the convexity parameter of $E^*$ cannot be larger than $\beta$.}\qed
\end{remark}

The estimation of the convexity parameter of $E^*$ in Theorem \ref{th:sc} is rather pessimistic when the number of alternatives increases. The dependence of $\beta$ on $n$ becomes explicit in case of independent and identically distributed random utility shocks.

\begin{corollary}[Strong convexity of $E^*$ for IID utility shocks]
\label{cor:iid}  
   Let the random utility shocks $\epsilon^{(i)}$, $i=1,\ldots,n$, be independent and identically distributed with the common probability density function $f$ having a mode $\bar z \in \R$.
   Then, the corresponding convex conjugate $E^*$ is $\beta$-strongly convex with respect to the $\|\cdot\|_1$ norm, where the convexity parameter is given by
\[
  \beta = \frac{1}{\displaystyle 2n(n-1) f\left(\bar z\right)}.
\]
\proof 
We estimate the probability distribution function of $\epsilon^{(j)}-\epsilon^{(i)}$:
\[
\begin{array}{rcl}
 \displaystyle g_{i,j} \left(z\right) &=&
  \displaystyle \int_{-\infty}^{\infty} \int_{-\infty}^{\infty} f_\epsilon \left(y^{(-i,j)} +x^{(i)}, z +x^{(i)},x^{(i)}\right) \diff x^{(i)} \diff y^{(-i,j)} \\ \\
  &=&
  \displaystyle
  \int_{-\infty}^{\infty} \int_{-\infty}^{\infty} \prod_{k \not = i,j} f\left(y_k +x^{(i)}\right)  f\left(z+x^{(i)}\right) f\left(x^{(i)}\right) \diff x^{(i)} \diff y^{(-i,j)} \\ \\
  &=&
  \displaystyle
  \int_{-\infty}^{\infty} f\left(z+x^{(i)}\right) f\left(x^{(i)}\right) \prod_{k \not = i,j}\int_{-\infty}^{\infty}  f\left(y_k +x^{(i)}\right) \diff y_{k}   \diff x^{(i)}.
  \end{array}
\]
For all $k\not =i,j$, and $x^{(i)} \in \R$ we have:
\[
  \int_{-\infty}^{\infty}  f\left(y_k +x^{(i)}\right) \diff y_{k} = 1.
\]
Moreover, $f\left(z+x^{(i)}\right) \leq f\left(\bar z \right)$ for all $x^{(i)} \in \R$. Altogether, we get
\[
 \displaystyle g_{i,j} \left(z\right) \leq f\left(\bar z \right).
\]
The application of Theorem \ref{th:sc} yields the assertion.
\qed
\end{corollary}

In what follows we concentrate on some special distributions of random utility shocks widely used in the discrete choice literature. For them we obtain better estimations of the convexity parameter of $E^*$, in particular, not dependent on the number of alternatives. First, we consider the class of generalized extreme value models \cite{mcfadden:1978}. The vector $\epsilon= \left(\epsilon^{(1)}, \ldots, \epsilon^{(n)}\right)^T$ of random utility shocks defines a generalized extreme value model (GEV) if it follows the joint distribution given by the probability density function
\[
  f_\epsilon\left(y^{(1)}, \ldots, y^{(n)}\right) = \frac{\partial^n \exp\left(-G\left(e^{-y^{(1)}},\ldots, e^{-y^{(n)}}\right)\right)}{\partial y^{(1)} \cdots \partial y^{(n)}},
\]
where the generating function $G:\R^n_+ \rightarrow \R_+$ has the following properties:
\begin{itemize}
\item[(G1)] $G$ is homogeneous of degree $\nicefrac{1}{\mu} > 0$.
\item[(G1)] $G\left(x^{(1)}, \ldots, x^{(i)}, \ldots, x^{(n)}\right) \rightarrow \infty$ as $x^{(i)} \rightarrow \infty$, $i=1, \ldots,n$.
\item[(G3)] %any subset of indices $\left\{ i_1, \ldots, i_k\right\} \subset \{1,\ldots,n\}$ 
For the partial derivatives of $G$ with respect to $k$ distinct variables it holds: 
\[
\frac{\partial^k G\left(x^{(1)},\ldots, x^{(n)}\right)}{\partial x^{\left(i_1\right)} \cdots \partial x^{\left(i_k\right)}} \geq 0 \mbox{ if } k \mbox{ is odd}, \quad 
\frac{\partial^k G\left(x^{(1)},\ldots, x^{(n)}\right)}{\partial x^{\left(i_1\right)} \cdots \partial x^{\left(i_k\right)}} \leq 0 \mbox{ if } k \mbox{ is even}.
\]
\end{itemize}
It is well known from \cite{mcfadden:1978} that the surplus function for GEV is
\[
   E(u) = \mu \ln G\left(e^{u}\right) + \mu \gamma,
\]
where $\gamma$ is Euler's constant and we set $e^u=\left(e^{u^{(1)}}, \ldots, e^{u^{(n)}}\right)^T$ for the sake of brevity.
The choice probability of the $i$-th alternative is given by the $i$-th partial derivative of the GEV surplus function $E$:
\[
  \P \left( u^{(i)} + \epsilon^{(i)} = \max_{1 \leq i \leq n} u^{(i)} + \epsilon^{(i)}\right)=\frac{\partial E(u)}{\partial u^{(i)}} = \mu \frac{\partial G\left(e^{u}\right)}{\partial x^{(i)}}\cdot \frac{e^{u^{(i)}}}{G\left(e^{u}\right)}.
\]

We state a sufficient condition for the strong convexity of $E^*$ in terms of the generating function $G$. 

\begin{theorem}[Strong convexity of $E^*$ for GEV]
\label{th:gev}
   Let a generating function $G$ for GEV satisfy the following inequality for all $x=\left(x^{(1)}, \ldots, x^{(n)}\right)^T\in \R^n_+$:
\[
   \sum_{i=1}^{n} \frac{\partial^2 G(x)}{\partial x^{{(i)}2}} \cdot x^{{(i)}2} \leq M \cdot G(x)
\]
with some constant $M \in \R$. Then, the corresponding convex conjugate $E^*$ is $\beta$-strongly convex with respect to the $\|\cdot\|_1$ norm, where the convexity parameter is given by
\[
  \beta = \frac{1}{2(\mu M +1)-\nicefrac{1}{\mu}}.
\]
\proof 
For the second order partial derivative of $E$ with respect to $u^{(i)}$ it holds:
\[
\begin{array}{rcl}
  \displaystyle \frac{\partial^2 E(u)}{\partial u^{{(i)}2}} &=& \displaystyle \mu \left( \frac{\partial G\left(e^{u}\right)}{\partial x^{(i)}}\cdot \frac{e^{u^{(i)}}}{G\left(e^{u}\right)}  \left(1-\frac{\partial G\left(e^{u}\right)}{\partial x^{(i)}} \cdot\frac{e^{u^{(i)}}}{G\left(e^{u}\right)} \right)
   + 
    \frac{\partial G^2\left(e^{u}\right)}{\partial x^{{(i)}2}}\cdot\frac{\left(e^{u^{(i)}}\right)^2}{G\left(e^{u}\right)} \right) \\ \\
    &=& \displaystyle  \frac{1}{\mu}\frac{\partial E(u)}{\partial u^{(i)}} \left( 1 - \frac{\partial E(u)}{\partial u^{(i)}}\right) + \left(1-\frac{1}{\mu}\right)\frac{\partial E(u)}{\partial u^{(i)}}
   + \mu
    \frac{\partial G^2\left(e^{u}\right)}{\partial x^{{(i)}2}}\cdot\frac{\left(e^{u^{(i)}}\right)^2}{G\left(e^{u}\right)}.
    \end{array}
\] 
Analogously, we obtain the mixed partial derivative of $E$ for $j\not = i$:
\[
\begin{array}{rcl}
  \displaystyle 
   \frac{\partial^2 E(u)}{\partial u^{(i)} \partial u^{(j)}} &=& \displaystyle 
   \mu \left( -\frac{\partial G\left(e^{u}\right)}{\partial x^{(i)}} \frac{e^{u^{(i)}}}{G\left(e^{u}\right)} \cdot  \frac{\partial G\left(e^{u}\right)}{\partial x^{(j)}} \frac{e^{u^{(j)}}}{G\left(e^{u}\right)} 
   + 
    \frac{\partial G^2\left(e^{u}\right)}{\partial x^{(i)} \partial x^{(j)}}\cdot \frac{e^{u^{(i)}} e^{u^{(j)}}}{G\left(e^{u}\right)} \right) \\ \\
     &=&\displaystyle -\frac{1}{\mu}\frac{\partial E(u)}{\partial u^{(i)}} \cdot  \frac{\partial E(u)}{\partial u^{(j)}} 
   + \mu
    \frac{\partial G^2\left(e^{u}\right)}{\partial x^{(i)} \partial x^{(j)}}\cdot\frac{e^{u^{(i)}} e^{u^{(j)}}}{G\left(e^{u}\right)}.
\end{array}
\]
Equivalently, we have in matrix form:
\[
  \nabla^2 E(u) = \frac{1}{\mu} R(u) + S(u),  
\]
where
\[
 \begin{array}{rcl} 
 R(u) &= & \displaystyle \diag{\nabla E(u)} - \nabla E(u) \cdot \nabla^T E(u), \\ \\
 S(u) &=& \displaystyle \left(1-\frac{1}{\mu}\right) \diag{\nabla E(u)} + \mu \frac{\diag{e^u} \cdot \nabla^2  G\left(e^u\right)\cdot \diag{e^u}}{G\left(e^u\right)}.
\end{array}
\]
Since $\nabla E(u) \in \Delta$, we may apply Corollary \ref{cor:pp} to derive that
\[
 R(u) \in \mathcal{A}, \quad \|R(u)\|_{\infty,1} \leq 1.
\]
In view of $\nabla^2E(u), R(u) \in \mathcal{A}$, the matrix $S(u)$ fulfills (A2):
\[
   S(u) \cdot e = \underbrace{\nabla^2 E(u)\cdot e}_{=0} - \frac{1}{\mu} \underbrace{R(u)\cdot e}_{=0} = 0.   
\]
The off-diagonal entries of $S(u)$ are nonpositive due to (G3), hence, $S(u)$ fulfills also (A1). We apply Lemma \ref{lem:tr} to estimate the $\|\cdot\|_{\infty,1}$ norm of $S(u)\in \mathcal{A}$:
\[
\begin{array}{rcl}
    \displaystyle  
  \|S(u)\|_{\infty,1} &\leq& 2 \tr{S(u)} \\ \\ &=& \displaystyle 2 \left(\left(1-\frac{1}{\mu}\right) \sum_{i=1}^{n} \frac{\partial E(u)}{\partial u^{(i)}} + \mu \sum_{i=1}^{n}  \frac{\partial G^2\left(e^{u}\right)}{\partial x^{{(i)}2}}\cdot\frac{\left(e^{u^{(i)}}\right)^2}{G\left(e^{u}\right)}\right) \\ \\
  &\leq& \displaystyle 2 \left(\left(1-\frac{1}{\mu}\right) + \mu M \right).
\end{array}
\]
Overall, we have:
\[
   \|\nabla^2 E(u)\|_{\infty,1} \leq \frac{1}{\mu} \|R(u)\|_{\infty,1} + \|S(u)\|_{\infty,1} \leq \frac{1}{\mu} + 2 \left(\left(1-\frac{1}{\mu}\right) + \mu M \right)=\frac{1}{\beta}.
\]
\end{theorem}
The latter provides, due to Lemma \ref{lem:c2}, that the surplus function $E$ is $\frac{1}{\beta}$-strongly smooth with respect to the maximum norm $\|\cdot\|_\infty$. Finally, we apply Lemma \ref{lem:sd} to conclude that the convex conjugate $E^*$ is $\beta$-strongly convex with respect to the $\|\cdot\|_1$ norm.\qed

Now, we consider generalized nested logit models (GNL) introduced in \cite{wen:2001}. GNL is a particular class of GEV models defined by the generating function
\[
   G(x)= \sum_{\ell \in L} \left( \sum_{i =1}^{n} \left(\sigma_{i\ell}\cdot x^{(i)}\right)^{\nicefrac{1}{\mu_\ell}} \right)^{\nicefrac{\mu_\ell}{\mu}}.
\]
Here, $L$ is a generic set of nests. The parameters $\sigma_{i\ell} \geq 0$ denote the shares of the $i$-th alternative with which it is attached to the $\ell$-th nest. For any fixed $i \in \{1, \ldots,n\}$ they sum up to one:
\[
  \sum_{\ell \in L} \sigma_{i\ell}  = 1,
\]
and $\sigma_{i\ell}=0$ means that the $\ell$-th nest does not contain the $i$-th alternative. Hence, the set of alternatives within the $\ell$-th nest is
\[
  N_\ell = \left\{i \,|\, \sigma_{i\ell} >0\right\}.
\]
The nest parameters $\mu_\ell > 0$ describe the variance of the random errors while choosing alternatives within the $\ell$-th nest. Analogously, $\mu >0$ describes the variance of the random errors while choosing among the nests. For the function $G$ to fulfill (G1)-(G3) we require:
\[
   \mu_\ell \leq \mu \quad \mbox{for all } \ell \in L. 
\]
The underlying choice process can be viewed to comprise two stages:
\begin{itemize}
    \item[(1)] the probability of choosing the $\ell$-th nest is
\[
     q_\ell = \frac{e^{\nicefrac{v_\ell}{\mu}}}{\displaystyle
     \sum_{\ell \in L}e^{\nicefrac{v_\ell}{\mu}}},
\]
where
\[
v_\ell = \mu_\ell \ln \left( \sum_{i =1}^{n} \left(\sigma_{i\ell} \cdot e^{u^{(i)}}\right)^{\nicefrac{1}{\mu_\ell}} \right)
\]
stands for the utility attached to the $\ell$-th nest;
 \item[(2)] the probability of choosing the $i$-th alternative within the $\ell$-th nest is
\[
     p_{i\ell} = \frac{\left(\sigma_{i\ell} \cdot e^{u^{(i)}}\right)^{\nicefrac{1}{\mu_\ell}}}{\displaystyle
     \sum_{i=1}^{n} \left(\sigma_{i\ell} \cdot e^{u^{(i)}}\right)^{\nicefrac{1}{\mu_\ell}}}.
\]
\end{itemize}
Overall, the choice probability of the $i$-th alternative according to GNL amounts to
\[
  \P \left( u^{(i)} + \epsilon^{(i)} = \max_{1 \leq i \leq n} u^{(i)} + \epsilon^{(i)}\right)= \mu \frac{\partial G\left(e^{u}\right)}{\partial x^{(i)}}\cdot\frac{e^{u^{(i)}}}{G\left(e^{u}\right)}= \sum_{\ell \in L} q_\ell \cdot p_{i\ell}.
\]

\begin{corollary}[Strong convexity of $E^*$ for GNL]
\label{cor:gnl}
   For GNL the corresponding convex conjugate $E^*$ is $\beta$-strongly convex with respect to the $\|\cdot\|_1$ norm, where the convexity parameter is given by
\[
  \beta = \frac{1}{\frac{2}{\displaystyle \min_{\ell \in L} \mu_\ell}-\nicefrac{1}{\mu}}.
\]
\proof Let us estimate the constant $M$ from Theorem \ref{th:gev}. We have:
\[
  \frac{\partial G\left(x\right)}{\partial x^{(i)}}= \frac{1}{\mu}\sum_{\ell \in L} \left( \sum_{i =1}^{n} \left(\sigma_{i\ell}\cdot x^{(i)}\right)^{\nicefrac{1}{\mu_\ell}} \right)^{\nicefrac{\mu_\ell}{\mu}-1}
   \left(\sigma_{i\ell}\cdot x^{(i)}\right)^{\nicefrac{1}{\mu_\ell}-1} \cdot \sigma_{i\ell},
\]
and, further:
\[
 \begin{array}{rcl}
 \displaystyle \frac{\partial^2 G(x)}{\partial x^{{(i)}2}} &=& \displaystyle \frac{1}{\mu}\sum_{\ell \in L} 
 \frac{1}{\mu_\ell}\left(\frac{\mu_\ell}{\mu}-1\right)
 \left( \sum_{i =1}^{n} \left(\sigma_{i\ell}\cdot x^{(i)}\right)^{\nicefrac{1}{\mu_\ell}} \right)^{\nicefrac{\mu_\ell}{\mu}-2}
   \left(\left(\sigma_{i\ell}\cdot x^{(i)}\right)^{\nicefrac{1}{\mu_\ell}-1} \cdot\sigma_{i\ell}\right)^2 \\ \\
  &+& \displaystyle \frac{1}{\mu}\sum_{\ell \in L} \left(\frac{1}{\mu_\ell}-1 \right) \left( \sum_{i =1}^{n} \left(\sigma_{i\ell}\cdot x^{(i)}\right)^{\nicefrac{1}{\mu_\ell}} \right)^{\nicefrac{\mu_\ell}{\mu}-1}
  \left(\sigma_{i\ell}\cdot x^{(i)}\right)^{\nicefrac{1}{\mu_\ell}-2}  \cdot\sigma_{i\ell}^2.
 \end{array}
\]
Due to $\mu_\ell \leq \mu$, $\ell \in L$, we get:
\[
  \frac{\partial^2 G(x)}{\partial x^{{(i)}2}} \leq \frac{1}{\mu}\sum_{\ell \in L} \left(\frac{1}{\mu_\ell}-1 \right) \left( \sum_{i =1}^{n} \left(\sigma_{i\ell}\cdot x^{(i)}\right)^{\nicefrac{1}{\mu_\ell}} \right)^{\nicefrac{\mu_\ell}{\mu}-1}
   \left(\sigma_{i\ell}\cdot x^{(i)}\right)^{\nicefrac{1}{\mu_\ell}-2} \cdot\sigma_{i\ell}^2.
\]
We multiply these terms by $x^{{(i)}2}$ and sum up over $i=1, \ldots, n$:
\[
\begin{array}{rcl}
 \displaystyle 
   \sum_{i=1}^{n} \frac{\partial^2 G(x)}{\partial x^{{(i)}2}} \cdot x^{{(i)}2} &\leq& \displaystyle
   \frac{1}{\mu}\sum_{\ell \in L} \left(\frac{1}{\mu_\ell}-1 \right) \left( \sum_{i =1}^{n} \left(\sigma_{i\ell}\cdot x^{(i)}\right)^{\nicefrac{1}{\mu_\ell}} \right)^{\nicefrac{\mu_\ell}{\mu}-1} \\ \\ &&\displaystyle \cdot 
   \sum_{i=1}^{n}\left(\sigma_{i\ell}\cdot x^{(i)}\right)^{\nicefrac{1}{\mu_\ell}-2} \cdot  \sigma_{i\ell}^2 \cdot x^{{(i)}2} \\ \\
  &=& \displaystyle \frac{1}{\mu}\sum_{\ell \in L} \left(\frac{1}{\mu_\ell}-1 \right) \left( \sum_{i =1}^{n} \left(\sigma_{i\ell}\cdot x^{(i)}\right)^{\nicefrac{1}{\mu_\ell}} \right)^{\nicefrac{\mu_\ell}{\mu}}\\ \\
   &\leq& \displaystyle\frac{1}{\mu} \max_{\ell \in L} \left( \frac{1}{\mu_\ell} -1\right) \cdot G(x)= \frac{1}{\mu}  \left( \frac{1}{\displaystyle \min_{\ell \in L}\mu_\ell} -1\right) \cdot G(x).
 \end{array}
\]
Hence,  we may apply Theorem \ref{th:gev} with 
\[
  M = \frac{1}{\mu}  \left( \frac{1}{\displaystyle \min_{\ell \in L}\mu_\ell} -1\right).
\]
\qed
\end{corollary}

\begin{example}[Multinomial logit] \textup{Let in GNL there be just one nest, i.\,e. $L = \{1\}$, and $\mu_1 = \mu$. Then, the generating function
\[
   G(x)= \sum_{i =1}^{n} \left(x^{(i)}\right)^{\nicefrac{1}{\mu}}
\]
leads to the multinomial logit (ML). The corresponding surplus function is
\[
   E(u) = \mu \ln \sum_{i =1}^{n} e^{\nicefrac{u^{(i)}}{\mu}} + \mu \gamma,
\]
and the choice probabilities are
\[
   \P \left( u^{(i)} + \epsilon^{(i)} = \max_{1 \leq i \leq n} u^{(i)} + \epsilon^{(i)}\right)= \frac{e^{\nicefrac{u^{(i)}}{\mu}}}{\displaystyle
     \sum_{i=1}^{n}e^{\nicefrac{u^{(i)}}{\mu}}}, \quad i=1, \ldots, n.
\]
Note that ML can be deduced from the IID random utility shocks $\epsilon^{(i)}$, $i=1, \ldots,n$, each of them following the Gumbel distribution with zero mode and variance $\nicefrac{\mu \pi}{\sqrt{6}}$ \cite{depalma:1994}. For ML the convex conjugate of the surplus function can be explicitly given:
\[
   E^*(p) = \mu \sum_{i=1}^{n} p^{(i)} \ln p^{(i)} - \mu \gamma = \mu H(p) - \mu \gamma,
\]
where $H$ is the (negative) entropy. It is well known that $H$ is $1$-strongly convex with respect to the $\|\cdot\|_1$ norm due to the Pinsker inequality. Hence, $E^*$ is $\mu$-strongly convex with respect to the $\|\cdot\|_1$ norm. The same result also follows from Corollary \ref{cor:gnl} with the convexity parameter
\[
  \beta = \frac{1}{\frac{2}{\displaystyle \min_{\ell \in L} \mu_\ell}-\nicefrac{1}{\mu}}= \frac{1}{\nicefrac{2}{\mu}-\nicefrac{1}{\mu}}= \mu.
\]
} \qed
\end{example}

\begin{example}[Nested logit]
\textup{Let in GNL for every alternative $i \in \{1,\ldots,n\}$ there be a unique nest $\ell_i \in L$ with $\sigma_{i \ell_i}=1$, and $\mu=1$. Then, 
the nests $N_\ell=\left\{i \,|\, \ell_i=\ell\right\}$ are mutually exclusive, and 
the generating function
\[
   G(x)=\sum_{\ell \in L} \left( \sum_{i \in N_\ell} x^{{(i)}\nicefrac{1}{\mu_\ell}} \right)^{\mu_\ell}
\]
leads to the nested logit (NL). The corresponding surplus function is
\[
   E(u) = \mu \ln \sum_{\ell \in L} \left( \sum_{i \in N_\ell} e^{\nicefrac{u^{(i)}}{\mu_\ell}} \right)^{\mu_\ell} + \mu \gamma,
\]
and the choice probabilities for $i \in N_\ell$, $\ell \in L$ are
\[
   \P \left( u^{(i)} + \epsilon^{(i)} = \max_{1 \leq i \leq n} u^{(i)} + \epsilon^{(i)}\right)= 
   \frac{e^{\mu_\ell \ln \sum_{i \in N_\ell} e^{\nicefrac{u^{(i)}}{\mu_\ell}}}}{\displaystyle
     \sum_{\ell \in L} e^{\mu_\ell \ln \sum_{i \in N_\ell} e^{\nicefrac{u^{(i)}}{\mu_\ell}}}}\cdot   \frac{e^{\nicefrac{u^{(i)}}{\mu_\ell}}}{\displaystyle
     \sum_{i\in N_\ell} e^{\nicefrac{u^{(i)}}{\mu_\ell}}}.
\]
For NL the convex conjugate of the surplus function is explicitly given in \cite{fosgerau:2017}:
\[
   E^*(p) = \sum_{\ell \in L} \mu_\ell \sum_{i\in N_\ell} p^{(i)} \ln p^{(i)} 
   + \sum_{\ell \in L} \left(1-\mu_\ell\right) \left(\sum_{i\in N_\ell} p^{(i)} \right) \ln \left(\sum_{i\in N_\ell} p^{(i)} \right)
   - \mu\gamma.
\]
By examining the first part of this formula, it can be shown that $E^*$ is  $\displaystyle \left(\min_{\ell \in L} \mu_\ell\right)$-strongly convex with respect to the $\|\cdot\|_1$ norm. Corollary \ref{cor:gnl} provides the convexity parameter, which is of the same order:
\[
  \beta = \frac{1}{\frac{2}{\displaystyle \min_{\ell \in L} \mu_\ell}-\nicefrac{1}{\mu}}= \frac{1}{\frac{2}{\displaystyle \min_{\ell \in L} \mu_\ell}-1} > \frac{1}{2} \min_{\ell \in L} \mu_\ell.
\]
} \qed
\end{example}

\begin{example}[GNL's from literature] \textup{For other specifications of GNL -- except of ML and NL -- the explicit form of the convex conjugate $E^*$ is not known yet. In this case Corollary \ref{cor:gnl} can be applied in order to estimate the convexity parameter of $E^*$. 
\begin{itemize}
    \item[(i)] {\bf Ordered GEV} \cite{small:1987} is a GNL model with
\[
  L=\{ 1, \ldots, n+m\}, \quad \mu=1,
\]
\[
 \sigma_{i\ell}>0 \mbox{ for all } \ell \in \{i, \ldots,i+m\}, \quad
 \sigma_{i\ell}=0 \mbox{ for all } \ell \in L \backslash \{i, \ldots,i+m\}.
\]
There are $n+m$ overlapping nests $N_\ell=\left\{i \,|\, \ell-m \leq i \leq \ell\right\}$, and every alternative lies exactly in $m+1$ of them, namely $i \in N_\ell$ for $\ell=i, \ldots, i+m$. Then, 
the generating function is 
\[
   G(x)=\sum_{\ell=1}^{n+m} \left( \sum_{i \in N_\ell} \left(\sigma_{i\ell} x^{(i)}\right)^{\nicefrac{1}{\mu_\ell}} \right)^{\mu_\ell}.
\]
    \item[(ii)] {\bf Paired combinatorial logit} \cite{koppelman:2000} is a GNL model with
\[
  L=\{ (i,j)\in \{1,\ldots, n\} \,|\,  i \not = j\}, \quad \mu=1,
\]
\[
\sigma_{i\ell}=
\left\{ \begin{array}{cl}
     \displaystyle\frac{1}{2(n-1)} & \mbox{if } \ell = (i,j), (j,i) \mbox{ with } j \not =i,  \\
     0 & \mbox{else.} 
\end{array}
\right.
\]
There are $n^2-n$ nests corresponding to the pairs of alternatives, and every alternative lies in $2(n-1)$ of them. Then, 
the generating function is 
\[
   G(x)=\sum_{\ell=(i,j), i \not =j} \left( \left(\sigma_{i\ell} x^{(i)}\right)^{\nicefrac{1}{\mu_\ell}} + \left(\sigma_{j\ell} x^{(j)}\right)^{\nicefrac{1}{\mu_\ell}} \right)^{\mu_\ell}.
\]
 \item[(iii)] {\bf Principles of differentiation GEV} \cite{bresnahan:1997} is a GNL model with
\[
  L = \mathop{\dot{\bigcup}}_{d \in D} L_d, \quad \mu=1, \quad
  \mu_\ell = \mu_d \mbox{ for all } \ell \in L_d,
\]
\[  
 \sigma_{i\ell}= 
\left\{ \begin{array}{cl}
     \displaystyle \sigma_d & \mbox{if } i \in N_{\ell d} \mbox{ and } \ell \in L_d, \\
     0 & \mbox{else,} 
\end{array}
\right.
\]
where
\[
  \{1,\ldots,n\}=\mathop{\dot{\bigcup}}_{\ell \in L_d} N_{\ell d}.
\]
The set $D$ represents the dimensions of alternatives. Along the $d$-th dimension alternatives can be clustered into the disjoint nests $N_{\ell d}$, $\ell \in L_d$. The shares $\sigma_{i\ell}$ of all alternatives within these nests depend only on the dimension $d$. The parameters $\mu_\ell$ also coincide for all $\ell \in L_d$.
Then, the generating function is 
\[
   G(x)=\sum_{d \in D} \sigma_d \sum_{\ell \in L_d} \left( \sum_{i\in N_{\ell d}} \left(x^{(i)}\right)^{\nicefrac{1}{\mu_d}} \right)^{\mu_d}.
\]
\end{itemize}
Note that the convexity parameter of $E^*$ for GNL specifications (i)-(iii) depends only on the smallest of the nest parameters $\mu_\ell$, $\ell \in L$:
\[
  \beta = \frac{1}{\frac{2}{\displaystyle \min_{\ell \in L} \mu_\ell}-1}.
\]
\qed
}
\end{example}

Let us relate our results on the strong convexity of $E^*$ to those existing in the literature.

\begin{remark}[Cross moment model]
In CMM \cite{mishra:2012} the random vector $\epsilon \sim (0,\Sigma)$
follows a joint distribution with zero mean and a given covariance matrix $\Sigma$. Additionally, it maximizes the surplus function:
\[
Z(u) = \max_{\epsilon \sim (0,\Sigma)} \E_\epsilon \left(\max_{1 \leq i \leq n} u^{(i)} + \epsilon^{(i)} \right).
\]
Assuming that the covariance matrix $\Sigma$ is positive definite, the following dual representation of $Z$ has been derived in \cite{li:2019}:
\[
      Z(u) =  \max_{ p \in \Delta} \, \langle p, u \rangle + \mbox{tr} \left( \left(\Sigma^{\nicefrac{1}{2}}\left( \mbox{diag}(p) - p p^T \right) \Sigma^{\nicefrac{1}{2}}\right)^{\nicefrac{1}{2}} \right).
\]
%Here, the unique square root of positive semidefinite matrices is used. 
Moreover, the solution $p\in \Delta$ of the latter optimization problem provides the choice probabilities corresponding to the random error $\epsilon(u)\sim (0,\Sigma)$ that maximizes the surplus function above, i.\,e.
\[
   p^{(i)} = \P \left( u^{(i)} + \epsilon^{(i)}(u) = \max_{1 \leq i \leq n} u^{(i)} + \epsilon^{(i)}(u)\right), \quad i=1, \ldots, n.
\]
Hence, the convex conjugate of $Z$ is
\[
   Z^*(p) = - \mbox{tr} \left( \left(\Sigma^{\nicefrac{1}{2}}\left( \mbox{diag}(p) - p p^T \right) \Sigma^{\nicefrac{1}{2}}\right)^{\nicefrac{1}{2}} \right).
\]
In \cite{li:2019}, $Z^*$ is shown to be strongly convex on the simplex w.r.t. Euclidean norm (however, its convexity parameter is not given explicitly). We point out that $Z^*$ can be therefore used as a discrete choice prox-function on the simplex as well. \qed
\end{remark}

\subsubsection{Simplicity}
\label{ssec:sim}

In view of the discussion in Section \ref{ssec:sc}, we assume that the convex conjugate $E^*$ is strongly convex. 

\begin{theorem}[Simplicity]
\label{th:sim}
   The unique maximizer of the optimization problem 
\begin{equation}
    \label{eq:sim}
    E(u) =  \sup_{ p \in \Delta} \, \langle p, u \rangle - E^*(p)
\end{equation}
is given by the choice probabilities
\[
    p^{(i)} = \P \left( u^{(i)} + \epsilon^{(i)} = \max_{1 \leq i \leq n} u^{(i)} + \epsilon^{(i)}\right), \quad i=1, \ldots, n.  
\]
\proof Let a vector $u \in \R^n$ of deterministic utilities be given. The vector $p$ of corresponding probabilities lies in the relative interior of the simplex $\Delta$, i.\,e. $p \in \mbox{rint}(\Delta)$. This is due to the fact that the distribution of the random vector $\epsilon$ is fully supported on $\R^n$, see \cite{norets:2013}. Hence, we have:
\[
  p = \nabla E(u),
\]
or, equivalently, by the convex duality \cite{rockafellar:1970}:
\[
  u \in \partial E^*(p).
\]
The latter is the first order optimality condition for  the optimization problem (\ref{eq:sim}). Additionally, note that 
the dual representation of the surplus function is valid due to the Fenchel-Moreau theorem.
\qed
\end{theorem}

\begin{corollary}[Lower bound for $E^*$] 
\label{cor:lb}
The unique minimizer $p_0$ of the convex conjugate $E^*$ consists of the choice probabilities with respect to the zero-utility, i.\,e.
\[
    p^{(i)}_0 = \P \left( \epsilon^{(i)} = \max_{1 \leq i \leq n} \epsilon^{(i)}\right), \quad i=1, \ldots, n.  
\]
Moreover, it holds:
\[
    E^*\left(p\right) \geq 
    E^*\left(p_0\right) = -  \E_\epsilon \left(\max_{1 \leq i \leq n}\epsilon^{(i)} \right) \quad \mbox{for all } p \in \Delta.
\]
\proof
Setting $u=0$ in Theorem \ref{th:sim} yields the assertion.
\qed
\end{corollary}

\section{Economic application}
\label{sec:ea}

We apply discrete choice prox-functions for the natural adjustment of consumer's demand.

\subsection{Lancaster's approach to consumer theory}
\label{sec:lan}

We briefly describe the Lancaster's approach to consumer theory as presented in \cite{lancaster:1966}. For that, let $x \in \R^m_+$ denote the consumer's demand vector of $m$ goods. Every such demand vector $x$ generates the vector $z \in \R^n$ of $n$ qualities (sometimes called characteristics):
\[
     z = Q x,
\]
where $Q=\left(q_{ij}\right) \in \R^{n\times m}$ is a fixed quality matrix. Its entries $q_{ij}$ denote the amounts of the $i$-th quality while consuming one unit of the $j$-th good. Further, the consumer assigns utility $u(z)$ to the goods' qualities $z$, and tries to maximize it by adjusting the demand $x$. Hereby, the budget constraint need to be satisfied, i.\,e.
\[
   \langle \pi, x \rangle \leq w,
\]
where $\pi \in \R^m_{+}$ is a fixed vector of positive goods' prices, i.\,e. $\pi >0$, and $w >0$ is a fixed available budget. Overall, the Lancaster's approach to consumer theory consists in solving the following maximization problem:
\[  
    \max_{\scriptsize \begin{array}{c} x \geq 0 \end{array}} \, u(z) \quad \mbox{s.t.} \quad z = Q x, \quad \langle \pi, x \rangle \leq w.
\]
In what follows we focus on the Leontieff utility function by setting
\[ 
   u(z)=\min \left\{ \frac{z^{(1)}}{\sigma^{(1)}}, \ldots, \frac{z^{(n)}}{\sigma^{(n)}}\right\},
\]
where $\sigma=\left(\sigma^{(1)}, \ldots, \sigma^{(n)}\right)^T\in \R^n_{+}$ is a fixed vector of positive quality standards, i.\,e. $\sigma >0$.
%Since the Leontieff utility function $u(Qx)$ is homogeneous of degree one w.r.t. $x$, we may assume without loss of generality that the budget is normalized, i.\,e. $w=1$. 
Additionally,
we assume that there exists a feasible demand vector which delivers positive Leontieff utility. Otherwise, the Lancaster's consumption problem is trivially solved by driving the demand to zero. On the contrary, under the proposed assumption the whole budget $w$ will be spent at any optimal demand. Hence, the budget constraint can be taken tight without loss of generality. The primal Lancaster's consumption problem becomes
\[
   \mbox{P:}\quad \max_{\scriptsize \begin{array}{c} x \geq 0\\\langle \pi, x \rangle = w   \end{array}} U(x),
\]
with the concave objective function
\[
U(x)= u(Qx)=\min_{1\leq i \leq n} \frac{(Qx)^{(i)}}{\sigma^{(i)}}.
\]
The optimization problem (P) consists in adjusting demand $x$ by spending the budget $w$ in order to maximize the worst ratio of consumed qualities $(Qx)^{(i)}$ in relation to their desired standards $\sigma^{(i)}$ taking over $i=1, \ldots, n$.

Let us derive a dual optimization problem for (P). For that, we introduce dual variables 
\[
\lambda=\left(\lambda^{(1)}, \ldots, \lambda^{(n)}\right)^T \in \R^n, 
\]
which can be interpreted as internal prices of qualities. Due to the duality of linear programming, we have:
\[
\begin{array}{rcl}
   \displaystyle  
   \max_{\scriptsize \begin{array}{c} x \geq 0\\\langle \pi, x \rangle = w   \end{array}} U(x) &=& \displaystyle
   \max_{\scriptsize \begin{array}{c} x \geq 0\\\langle \pi, x \rangle = w   \end{array}}  \min_{1\leq i \leq n} \frac{(Qx)^{(i)}}{\sigma^{(i)}} \\ \\ &=&  \displaystyle 
%  = \max_{\scriptsize \begin{array}{c} x \geq 0\\\langle p, x \rangle = w   \end{array}} \, \min_{\scriptsize \begin{array}{c} \lambda \geq 0\\\langle e, \lambda \rangle = 1   \end{array}} \sum_{i=1}^{n}\frac{(Qx)_i}{\sigma_i} \lambda_i 
   \max_{\scriptsize \begin{array}{c} x \geq 0\\\langle \pi, x \rangle = w   \end{array}}  \min_{\scriptsize \begin{array}{c} \lambda \geq 0\\\langle \sigma, \lambda \rangle = 1   \end{array}} \langle Qx, \lambda \rangle \\ \\
   &=& \displaystyle
    \min_{\scriptsize \begin{array}{c} \lambda \geq 0\\\langle \sigma, \lambda \rangle = 1   \end{array}}   \max_{\scriptsize \begin{array}{c} x \geq 0\\\langle \pi, x \rangle = w   \end{array}} \left\langle x, Q^T\lambda \right\rangle \\ \\ &=& \displaystyle  
    \min_{\scriptsize \begin{array}{c} \lambda \geq 0\\\langle \sigma, \lambda \rangle = 1   \end{array}}  w \max_{1\leq j \leq m} \frac{\left(Q^T \lambda\right)^{(j)}}{\pi^{(j)}} = \displaystyle 
    \min_{\scriptsize \begin{array}{c} \lambda \geq 0\\\langle \sigma, \lambda \rangle = 1   \end{array}} \Phi(\lambda).
\end{array}
\]
The dual Lancaster's consumption problem becomes
\[
   \mbox{D:}\quad \min_{\scriptsize \begin{array}{c} \lambda \geq 0\\\langle \sigma, \lambda \rangle = 1   \end{array}}  \Phi(\lambda),
\]
with the convex objective function
\[
\Phi(\lambda)=w \max_{1\leq j \leq m} \frac{\left(Q^T \lambda\right)^{(j)}}{\pi^{(j)}}.
\]
The optimization problem (D) consists in adjusting internal prices of qualities $\lambda$ in order to minimize the best ratio of quality estimates  $\left(Q^T \lambda \right)^{(j)}$ of goods in relation to their market prices $\pi^{(j)}$ taking over $j=1, \ldots, m$.
Note that usually the number $n$ of goods  is considerably bigger than the number $m$ of their qualities. From a consumer it will be thus more plausible to expect that (D) is successively solved rather than (P).    

\subsection{Dual Averaging Scheme}
\label{sec:da}
We rewrite the dual optimization problem (D) by introducing new variables
\[
    p^{(i)} = \sigma^{(i)} \lambda^{(i)}, \quad i =1, \ldots, n.
\]
In terms of the variables $p=\left(p^{(1)},\ldots,p^{(n)} \right)^T$ we equivalently obtain the following auxiliary optimization problem on the simplex:
\[
   \mbox{A:}\quad \min_{\scriptsize \begin{array}{c} p \in \Delta   \end{array}}  \Psi\left(p\right),
\]
where the objective function is 
\[
\Psi(p)= \Phi\left(\frac{p}{\sigma}\right).
\]
For solving (A) we apply the dual averaging scheme from \cite{nesterov:2013}. For that, we use the family of prox-functions on the simplex discussed in Section \ref{sec:prox}:
\[
   d(p)=E^*(p) - E^*\left(p_0\right). 
\]
Note that $d$ is continuous on the simplex $\Delta$ (Theorem \ref{th:ct}), strongly convex with convexity parameter $\beta>0$ (Theorems \ref{th:sc} and \ref{th:gev}, Corollary \ref{cor:gnl}). The computation of its convex conjugate 
$
   d^*(u) = E(u) - E\left(0\right)
$
is simple (Theorem \ref{th:sim}). Additionally, Corollary \ref{cor:lb} provides us with the prox-center of $\Delta$:
\[
   p_0 = \mbox{arg} \min_{p \in \Delta} \, d(p).
\]
In view of $d\left(p_0\right)=0$, the dual averaging scheme can be initialized by $p_0$.
%with the parameters:
%\[
%   \lambda[t]=1, \quad \beta[t+1]=(t+1)\mu[t+1] \quad 
%   \mbox{for } t=0,1,\ldots
%\]

\[
% \tag{DA}\label{eq:method}
   \begin{tabular}{|c|}
      \hline \\
      \begin{tabular}{c}
         {\bf Dual Averaging Scheme for (A)} 
      \end{tabular}\\ \hline \\
      \begin{tabular}{ll}
            {\bf 1.} Compute $\nabla \Psi\left(p_k\right)$. \\ \\
            {\bf 2.} Set $\displaystyle s_{k+1}=\frac{1}{k+1}\sum_{\ell=0}^{k} \nabla \Psi\left(p_\ell\right)$. \\ \\
            {\bf 3.} Update 
            $\displaystyle p_{k+1} = \mbox{arg} \min_{p \in \Delta} \left\{ \left\langle s_{k+1}, p \right\rangle + \frac{d(p)}{\sqrt{k+1}}  \right\}$.\\
      \end{tabular}
      \\ \\ \hline
   \end{tabular}
\]
The convergence properties of the dual averaging scheme follow from \cite[Theorem 1]{nesterov:2013}. For $k\geq 0$ it holds:
\begin{equation}
    \label{eq:ineq}
  \delta_k\leq \frac{D}{\sqrt{k+1}}+\frac{M^2}{2 \beta} \cdot \frac{1}{k+1}\left(1+\sum_{\ell=1}^{k} \frac{1}{\sqrt{\ell}}\right),
\end{equation}
where
\[
  \delta_k = \max_{p \in \Delta} \, \frac{1}{k+1} \sum_{\ell=0}^{k} \left\langle \nabla \Psi\left(p_\ell\right), p_\ell - p \right\rangle, \quad 
  M = \max_{p \in \Delta} \left\|\nabla \Psi\left(p\right) \right\|_\infty, \quad D\geq \max_{p \in \Delta} \, d(p). 
\]

\subsection{Consumption cycle}
\label{sec:cyc}
Let us state the dual averaging scheme from Section \ref{sec:da} in terms of the primal and dual Lancaster consumption problem (P) and (D), respectively. 

\begin{itemize}
    \item[\bf 1.] We compute
    \[
       \nabla \Psi\left(p_k\right) = \nabla \Phi\left(\frac{p_k}{\sigma}\right) = \frac{\nabla \Phi \left(\lambda_k\right)}{\sigma}, 
    \]
    where $\lambda_k$ denotes the internal prices of qualities at the $k$-th iteration. It holds:
    \[
      \nabla \Phi \left(\lambda_k\right) =  
      \nabla \left( w \max_{1\leq j \leq m} \frac{\left(Q^T \lambda_k\right)^{(j)}}{\pi^{(j)}} \right) = w \frac{Q y_k}{\pi},
    \]
    where the sharing vector $y_k \in \Delta$ fulfills 
    \[
    y^{(j)}_k = 0 \quad \mbox{for } j \not \in J\left(\lambda_k\right)
    \]
    with the active index set
    \[
       J\left(\lambda_k\right) = \left\{ j \in \{1, \ldots, m\} \, \left| \,  \frac{\left(Q^T \lambda_k\right)^{(j)}}{\pi^{(j)}}=\max_{1\leq j \leq m} \frac{\left(Q^T \lambda_k\right)^{(j)}}{\pi^{(j)}} \right.\right\}.  
    \]
    In other words, $J\left(\lambda_k\right)$ contains goods with the best quality/price ratio estimated by means of internal prices $\lambda_k$.
    We set the demand at the $k$-th iteration as
    \[ 
       x_k = w\frac{y_k}{\pi}.
    \]
    Note that $x_k$ is feasible for (P), since
    \[
       \left\langle \pi,x_k\right\rangle = \left\langle \pi,w\frac{y_k}{\pi}\right\rangle = w \left\langle e,y_k\right\rangle =w.
    \]
    Moreover, the demand is concentrated on the goods with the best quality/price ratio, i.\,e. $x^{(j)}_k \not = 0$ if and only if $j \in J\left(\lambda_k\right)$. Overall, we obtain:
    \[
       \nabla \Psi\left(p_k\right) = \frac{Q x_k}{\sigma},
    \]
    which are the ratios of the consumed qualities $Q x_k$ in relation to their standards $\sigma$.
    \item[\bf 2.] We set
    \[
      s_{k+1}=\frac{1}{k+1}\sum_{\ell=0}^{k} \nabla \Psi\left(p_\ell\right) = \frac{1}{k+1}\sum_{\ell=0}^{k} \frac{Q x_\ell}{\sigma} 
      = \frac{Q \bar x_k}{\sigma}
    \]  
    with the average demand
    \[
       \bar x_k = \frac{1}{k+1}\sum_{\ell=0}^{k} x_\ell.
    \]
    Again, $s_{k+1}$ relates the average consumption $Q \bar x_k$ to the standards $\sigma$.
    \item[\bf 3.] We update
    \[
    p_{k+1} = \mbox{arg} \min_{p \in \Delta} \left\{ \left\langle s_{k+1}, p \right\rangle + \frac{d(p)}{\sqrt{k+1}} \right\}.
    \]
    Due to Theorem \ref{th:sim}, we equivalently obtain for $i=1, \ldots, n$:
    \[
      p_{k+1}^{(i)} =  \P \left(  s_{k+1}^{(i)} - \frac{\epsilon^{(i)}}{\sqrt{k+1}} = \min_{1 \leq i \leq n}   s_{k+1}^{(i)} - \frac{\epsilon^{(i)}}{\sqrt{k+1}}\right). 
    \]
    For the internal prices we have:
    \[
      \lambda_{k+1}^{(i)} =  \frac{1}{\sigma^{(i)}} \P \left( s_{k+1}^{(i)} - \frac{\epsilon^{(i)}}{\sqrt{k+1}} = \min_{1 \leq i \leq n}  s_{k+1}^{(i)} - \frac{\epsilon^{(i)}}{\sqrt{k+1}}\right).
    \]
%    Thus, $p_{k+1}^{(i)}$ has a sense of relative importance of quality $i$ for the consumer after $k$ iterations of the consumption process. 
    Thus, the internal price $\lambda_{k+1}^{(i)}$ of the $i$-th quality is proportional to the probability of detecting its average consumption $\left(Q \bar x_k\right)^{(i)}$ as the lowest one in comparison to the standard $\sigma^{(i)}$. Moreover, this detecting process is in accordance with additive random utility models allowing behavioral interpretations. E.\,g., the additive random errors $\frac{\epsilon^{(i)}}{\sqrt{k+1}}$ are diminishing with respect to the iteration number $k$, which accounts for the learning effect. 
    We further mention that $\lambda_{k+1}$ is feasible for (D), since
    \[
       \left\langle \sigma,\lambda_{k+1}\right\rangle = \left\langle e,p_{k+1}\right\rangle =1.
    \]
    So are the average internal prices 
\[
 \bar \lambda_k = \frac{1}{k+1}\sum_{\ell=0}^{k} \lambda_\ell.
\]
\end{itemize}

Let us examine the convergence properties of the proposed consumption cycle. 

\begin{theorem}[Consumption cycle]
The duality gap between (P) and (D) evaluated at the average demand and the average internal prices is closing at the optimal rate $\mbox{O}\left(\frac{1}{\sqrt{k+1}}\right)$. Namely, it holds for $k \geq 0$:
\[
0 \leq \Phi\left( \bar \lambda_k\right) - U\left(\bar x_k\right) \leq
\left(D+\frac{M^2}{\beta}\right) \frac{1}{ \sqrt{k+1}},
\]
where
\[
   M = w \max_{\scriptsize \begin{array}{c} 1 \leq i \leq n\\ 1 \leq j \leq m   \end{array}} \frac{\left|q_{i,j}\right|}{\sigma^{(i)} \cdot \pi^{(j)}}, \quad 
   D = \E_\epsilon \left(\max_{1 \leq i \leq n}\epsilon^{(i)} \right) - \min_{1 \leq i \leq n} \E_\epsilon \left(\epsilon^{(i)} \right).
\]
\end{theorem}
\proof
We compute the gap bound
\[
 \begin{array}{rcl}
 \displaystyle \delta_k &=& \displaystyle \max_{p \in \Delta} \, \frac{1}{k+1} \sum_{\ell=0}^{k} \left\langle \nabla \Psi\left(p_\ell\right), p_\ell - p \right\rangle =
  \frac{1}{k+1} \sum_{\ell=0}^{k} w \left\langle \frac{Q y_\ell}{\pi}, \lambda_\ell\right\rangle - \min_{p \in \Delta} \left\langle s_{k+1}, p \right\rangle \\ \\
  &=& \displaystyle  \frac{1}{k+1} \sum_{\ell=0}^{k} w \sum _{j \in J\left(\lambda_\ell \right)} y_\ell^{(j)} \cdot \frac{\left(Q^T \lambda_\ell\right)^{(j)}}{\pi^{(j)}} - \min_{1 \leq i \leq n} s_{k+1}^{(i)} \\ \\ &=&\displaystyle \frac{1}{k+1} \sum_{\ell=0}^{k} w \max_{1 \leq j \leq m} \frac{\left(Q^T \lambda_\ell\right)^{(j)}}{\pi^{(j)}}- \min_{1 \leq i \leq n} \frac{\left(Q \bar x_k\right)^{(i)}}{\sigma^{(i)}} = \frac{1}{k+1} \sum_{\ell=0}^{k} \Phi\left( \lambda_\ell\right) - U\left(\bar x_k\right).
\end{array}
 \]
The uniform bound for subgradients of $\Psi$ is given by
\[
M=\max_{p \in \Delta} \left\|\nabla \Psi\left(p\right) \right\|_\infty = \max_{\scriptsize \begin{array}{c} x \geq 0\\\langle \pi, x \rangle = w   \end{array}} \left\| \frac{Q x}{\sigma} \right\|_\infty =  w \max_{\scriptsize \begin{array}{c} 1 \leq i \leq n\\ 1 \leq j \leq m   \end{array}} \frac{\left|q_{i,j}\right|}{\sigma^{(i)} \cdot \pi^{(j)}}. 
\]
By applying Corollaries 
\ref{cor:ub} and \ref{cor:lb}, we estimate the constant $D$: 
\[ 
  \max_{p \in \Delta} \, d(p) = \max_{p \in \Delta} \,
  E^*(p) - E^*\left(p_0\right) \leq  - \min_{1 \leq i \leq n} \E_\epsilon \left(\epsilon^{(i)} \right) +\E_\epsilon \left(\max_{1 \leq i \leq n}\epsilon^{(i)} \right)= D.
\]
It is straightforward to see that
\[
  \frac{1}{k+1}\left(1+\sum_{\ell=1}^{k} \frac{1}{\sqrt{\ell}}\right) \leq \frac{2}{\sqrt{k+1}}.
\]
Finally, by using (\ref{eq:ineq}) we get:
\[
\begin{array}{rcl}
  \displaystyle \Phi\left( \bar \lambda_k\right) - U\left(\bar x_k\right) &\leq& \displaystyle \frac{1}{k+1} \sum_{\ell=0}^{k} \Phi\left( \lambda_\ell\right) - U\left(\bar x_k\right) = \delta_k \\ \\ &\leq& \displaystyle
\frac{D}{\sqrt{k+1}}+\frac{M^2}{2 \beta} \cdot \frac{1}{k+1}\left(1+\sum_{\ell=1}^{k} \frac{1}{\sqrt{\ell}}\right)=\left(D+\frac{M^2}{\beta}\right) \frac{1}{ \sqrt{k+1}}.
\end{array}
\]
\qed 
%\section*{Appendix}

\end{document}